\documentclass[12pt]{article}

    \usepackage{amsmath,amsthm, amsfonts,amssymb}

    \usepackage{graphicx}
    \usepackage{pstricks,pst-plot}
    \usepackage[vcentermath]{youngtab}

    \makeatletter
    \newfont{\footsc}{cmcsc10 at 8truept}
    \newfont{\footbf}{cmbx10 at 8truept}
    \newfont{\footrm}{cmr10 at 10truept}

    \makeatother
    \newtheorem{thm}{Theorem}[section]

    \newtheorem{cor}[thm]{Corollary}

    \theoremstyle{definition}
        \newtheorem{defn}[thm]{Definition}
        \newtheorem{exgr}[thm]{Example}

        \newtheorem{alg}[thm]{Algorithm}

    \theoremstyle{remark}

    \numberwithin{equation}{section}
    \numberwithin{figure}{section}

    \title{A Geometric Form for the Extended\\ Patience Sorting
    Algorithm}

    \author{Alexander Burstein\\
    \small Department of Mathematics\\[-0.8ex]
    \small Iowa State University\\[-0.8ex]
    \small Ames, IA 50011-2064, USA\\[-0.8ex]
    \small \texttt{burstein@math.iastate.edu}\\[1.6ex]
    Isaiah Lankham\thanks{The work of the second author was
    supported in part by the U.S. National Science Foundation
    under Grants DMS-0135345 and DMS-0304414.}\\
    \small Department of Mathematics\\[-0.8ex]
    \small University of California, Davis\\[-0.8ex]
    \small Davis, CA 95616-8633, USA\\[-0.8ex]
    \small \texttt{issy@math.ucdavis.edu}}

    \date{\small Submitted: June 1, 2005; Accepted: July 19, 2005 }

\begin{document}

    \maketitle

    \begin{abstract}

        Patience Sorting is a combinatorial algorithm that can be
        viewed as an iterated, non-recursive form of the Schensted
        Insertion Algorithm.  In recent work the authors extended
        Patience Sorting to a full bijection between the symmetric
        group and certain pairs of combinatorial objects (called
        \emph{pile configurations}) that are most naturally defined in
        terms of generalized permutation patterns and barred pattern
        avoidance.  This Extended Patience Sorting Algorithm is very
        similar to the Robinson-Schensted-Knuth (or RSK)
        Correspondence, which is itself built from repeated
        application of the Schensted Insertion Algorithm.

        In this work we introduce a geometric form for the Extended
        Patience Sorting Algorithm that is in some sense a natural
        dual algorithm to G. Viennot's celebrated Geometric RSK
        Algorithm.  Unlike Geometric RSK, though, the lattice paths
        coming from Patience Sorting are allowed to intersect.  We
        thus also give a characterization for the intersections of
        these lattice paths in terms of the pile configurations
        associated with a given permutation under the Extended
        Patience Sorting Algorithm.

        \bigskip\bigskip\bigskip\bigskip\bigskip

    \end{abstract}

    \section{Introduction}
    \label{sec:Introduction}

    The term \emph{Patience Sorting} was introduced in 1962 by C.L.
    Mallows \cite{refMallows1962, refMallows1963} as the name of a
    card sorting algorithm invented by A.S.C. Ross.  This algorithm
    works by first partitioning a shuffled deck of $n$ cards (which we
    take to be a permutation $\sigma \in \mathfrak{S}_{n}$) into
    sorted subsequences $r_{1}, r_{2}, \ldots, r_{m}$ called
    \emph{piles} and then gathering the cards up in order from the
    tops of these piles.  The procedure used in forming $r_{1}, r_{2},
    \ldots, r_{m}$ can be viewed as an iterated, non-recursive form of
    the Schensted Insertion Algorithm for interposing values into the
    rows of a Young tableau (see \cite{refAD1999} and
    \cite{refBLFPSAC05}).  Given $\sigma \in \mathfrak{S}_{n}$, we
    call this resulting collection of piles (given as part of the more
    general Algorithm~\ref{alg:ExtendedPSalgorithm} below) the
    \emph{pile configuration} corresponding to $\sigma$ and denote it
    by $R(\sigma)$.

    Given a pile configuration $R$, one forms its \emph{reverse
    patience word} $RPW(R)$ by listing the piles in $R$ ``from bottom
    to top, left to right'' (i.e., by reversing the so-called
    ``far-eastern reading'') as illustrated in
    Example~\ref{eg:PileConfigurationExample} below.  In recent work
    \cite{refBLFPSAC05} the authors used G. Viennot's (northeast)
    shadow diagram construction (defined in \cite{refViennot1977} and
    summarized in Section~\ref{sec:GeometricRSK:NEshadowDiagrams}
    below) to characterize these words in terms of the following
    pattern avoidance condition: Given $\sigma \in \mathfrak{S}_{n}$,
    each instance of the generalized permutation pattern
    $2\textrm{-}31$ in $RPW(R(\sigma))$ must be contained within an
    instance of the pattern $3\textrm{-}1\textrm{-}42$.  We call this
    restricted form of the generalized permutation pattern
    $2\textrm{-}31$ a \emph{(generalized) barred permutation pattern}
    and denote it by $3\textrm{-}\bar{1}\textrm{-}42$.  This
    notational convention is due to J. West, et al., and first
    appeared in the study of two-stack sortable permutations
    \cite{refDGG1998, refDGW1996, refWest1990}.  As usual, we denote
    the set of permutations $\sigma \in \mathfrak{S}_{n}$ that avoid
    the pattern $3\textrm{-}\bar{1}\textrm{-}42$ by
    $S_{n}(3\textrm{-}\bar{1}\textrm{-}42)$.  (See B\'{o}na
    \cite{refBona2004} for a review of permutation patterns in
    general.)

    \begin{exgr}\label{eg:PileConfigurationExample}
        Let $\sigma = 6 4 5 1 8 7 2 3 \in \mathfrak{S}_{8}$.  Then,
        using a simplified form of
        Algorithm~\ref{alg:ExtendedPSalgorithm} below, $\sigma$ has
        the pile configuration $R(\sigma) = \{\{6 > 4 > 1\}, \{5 >
        2\}, \{8 > 7 > 3\}\}$, which is visually represented as the
        following array of numbers:\smallskip

        \begin{center}
            \begin{tabular}{cc}
                \begin{minipage}[c]{36pt}
                    $R(\sigma) \ \ = $
                \end{minipage}
                 &
                \begin{minipage}[c]{100pt}
                    \begin{tabular}{l l l}
                        1 &    & 3 \\
                        4 & 2 & 7 \\
                        6 & 5 & 8
                    \end{tabular}
                \end{minipage}
            \end{tabular}
        \end{center}\smallskip

        \noindent Furthermore, $RPW(R(6 4 5 1 8 7 2 3)) = 64152873 \in
        S_{8}(3\textrm{-}\bar{1}\textrm{-}42)$.

    \end{exgr}

    In \cite{refBLFPSAC05} the authors also extended the process of
    forming piles under Patience Sorting so that it essentially
    becomes a full non-recursive analog of the famous
    Robinson-Schensted-Knuth (or RSK) Correspondence.  As with RSK,
    this Extended Patience Sorting Algorithm
    (Algorithm~\ref{alg:ExtendedPSalgorithm} below) takes a simple
    idea (that of placing cards into piles) and uses it to build a
    bijection between elements of the symmetric group
    $\mathfrak{S}_{n}$ and certain pairs of combinatorial objects.  In
    the case of RSK, one uses the Schensted Insertion Algorithm to
    build a bijection with pairs of standard Young tableaux having the
    same shape (a partition $\lambda$ of $n$, denoted $\lambda \vdash
    n$; see \cite{refSagan2000}).  However, in the case of Patience
    Sorting, one achieves a bijection between permutations and
    somewhat more restricted pairs of pile configurations.  In
    particular, these pairs must not only have the same shape (a
    composition $\gamma$ of $n$, denoted $\gamma \
    \reflectbox{$\multimap$} \ n$) but their reverse patience words
    must also simultaneously avoid containing certain generalized
    permutation patterns in the same positions (see
    \cite{refBLFPSAC05} for more details).  This restriction can also
    be understood geometrically using Viennot's (northeast) shadow
    diagram construction for the permutation implicitly defined by a
    pair of pile configurations (as is discussed in
    \cite{refBLFPSAC05}), but a full geometric characterization
    requires the dual southwest shadow diagram construction used to
    define Geometric Patience Sorting in
    Section~\ref{sec:GeometricPS}.\\

    Viennot introduced the shadow diagram of a permutation in the
    context of studying the Sch\"{u}tzenberger Symmetry Property for
    RSK (first proven using a direct combinatorial argument in
    \cite{refSchutzenberger1963}).  Specifically, one can use
    recursively defined shadow diagrams to construct the RSK
    Correspondence completely geometrically via a sequence of
    recursively defined collections of non-intersecting lattice paths
    (with such collections called ``shadow diagrams'').  Then, using a
    particular labelling of the constituent ``shadow lines'' in each
    shadow diagram, one recovers successive rows in the usual RSK
    insertion and recording tableaux.  The Sch\"{u}tzenberger Symmetry
    Property for RSK then immediately follows since reflecting these
    shadow diagrams through the line ``$y = x$'' both inverts the
    permutation and exactly interchanges the labellings on the shadow
    lines that yield the rows in these tableaux.

    We review Viennot's Geometric RSK Algorithm in
    Section~\ref{sec:GeometricRSK} below.  Then, in
    Section~\ref{sec:GeometricPS}, we define a natural dual to
    Viennot's construction that similarly produces a geometric
    characterization of the Extended Patience Sorting Algorithm.  As
    with RSK, the Sch\"{u}tzenberger Symmetry Property is then
    immediate for Extended Patience Sorting.  Unlike Geometric RSK,
    though, the lattice paths formed under Geometric Patience Sorting
    are allowed to intersect.  Thus, having defined these two
    algorithms, we classify in Section~\ref{sec:PSvsRSK} the types of
    intersections that can occur under Geometric Patience Sorting and
    then characterize when they occur in terms of the pile
    configurations associated to a given permutation under Extended
    Patience Sorting (Algorithm~\ref{alg:ExtendedPSalgorithm}
    below).\\

    \noindent We close this introduction by stating the Extending
    Patience Sorting Algorithm and giving a complete example.

    \begin{alg}[Extended Patience Sorting
    Algorithm]\label{alg:ExtendedPSalgorithm} Given a shuffled
    deck of cards $\sigma = c_{1} c_{2} \cdots c_{n}$, inductively
    build \emph{insertion piles} $R = R(\sigma) = \{r_{1}, r_{2},
    \ldots, r_{m}\}$ and \emph{recording piles} $S = S(\sigma) =
    \{s_{1}, s_{2}, \ldots, s_{m}\}$ as follows:

        \begin{itemize}
            \item Place the first card $c_{1}$ from the deck into
            a pile $r_{1}$ by itself, and set $s_{1} = \{1\}$.

            \item For each remaining card $c_{i}$ ($i = 2, \ldots,
            n$), consider the cards $d_{1}, d_{2}, \ldots, d_{k}$ atop
            the piles $r_{1}, r_{2}, \ldots, r_{k}$ that have already
            been formed.

                \begin{itemize}
                    \item If $c_{i} > \max\{d_{1}, d_{2}, \ldots,
                    d_{k}\}$, then put $c_{i}$ into a new pile
                    $r_{k+1}$ by itself and set $s_{k+1} =
                    \{i\}$.\medskip

                    \item Otherwise, find the left-most card
                    $d_{j}$ that is larger than $c_{i}$ and put
                    the card $c_{i}$ atop pile $r_{j}$ while
                    simultaneously putting $i$ at the bottom of
                    pile $s_{j}$.\\
                \end{itemize}
        \end{itemize}
    \end{alg}

    \begin{exgr}\label{eg:ExtendedPSexample}
        Let $\sigma = 6 4 5 1 8 7 2 3 \in \mathfrak{S}_{8}$.  Then
        according to Algorithm~\ref{alg:ExtendedPSalgorithm} we
        simultaneously form the following pile configurations:\\

        \begin{tabular}{l p{56pt} p{72pt} l p{56pt} l}
            \begin{minipage}[c]{48pt}$\phantom{foo}$\end{minipage}
            &
            \begin{minipage}[c]{48pt}insertion piles\end{minipage}
            &
            \begin{minipage}[c]{48pt}recording piles\end{minipage}
            &
            \begin{minipage}[c]{48pt}$\phantom{foo}$\end{minipage}
            &
            \begin{minipage}[c]{48pt}insertion piles\end{minipage}
            &
            \begin{minipage}[c]{48pt}recording piles\end{minipage}
         \end{tabular}\\

        \begin{tabular}{l p{56pt} p{72pt} l p{56pt} l}
            \begin{minipage}[c]{48pt}
                Form a new pile with \textbf{6}:
            \end{minipage}
            &
            \begin{tabular}{l l l}
                                  & & \\
                                  & & \\
                   \textbf{6} & &
            \end{tabular}
            &
            \begin{tabular}{l l l}
                                 & & \\
                                 & & \\
                  \textbf{1} & &
            \end{tabular}
            &
            \begin{minipage}[c]{48pt}
                Then play the \textbf{4} on it:
            \end{minipage}
            &
            \begin{tabular}{l l l}
                                  & & \\
                   \textbf{4} & & \\
                              6 & &
            \end{tabular}
            &
            \begin{tabular}{l l l}
                               & & \\
                           1 & & \\
                \textbf{2} & &
            \end{tabular}
        \end{tabular}\\ \\

        \begin{tabular}{l p{56pt} p{72pt} l p{56pt} l}
            \begin{minipage}[c]{48pt}
                Form a new pile with \textbf{5}:
            \end{minipage}
            &
            \begin{tabular}{l l l}
                   & & \\
                4 & & \\
                6 & \textbf{5} &
            \end{tabular}
            &
            \begin{tabular}{l l l}
                  & & \\
               1 & & \\
               2 & \textbf{3} &
            \end{tabular}
            &
            \begin{minipage}[c]{48pt}
                Add the \textbf{1} to left pile:
            \end{minipage}
            &
            \begin{tabular}{l l l}
                \textbf{1} & & \\
                           4 & & \\
                           6 & 5 &
            \end{tabular}
            &
            \begin{tabular}{l l l}
                           1 & & \\
                           2 & & \\
                \textbf{4} & 3 &
            \end{tabular}
        \end{tabular}\\ \\

        \begin{tabular}{l p{56pt} p{72pt} l p{56pt} l}
            \begin{minipage}[c]{48pt}
                Form a new pile with \textbf{8}:
            \end{minipage}
            &
            \begin{tabular}{l l l}
                1 & & \\
                4 & & \\
                6 & 5 & \textbf{8}
            \end{tabular}
            &
            \begin{tabular}{l l l}
               1 & & \\
               2 & & \\
               4 & 3 & \textbf{5}
            \end{tabular}
            &
            \begin{minipage}[c]{48pt}
                Then play the \textbf{7} on it:
            \end{minipage}
            &
            \begin{tabular}{l l l}
                1 & & \\
                4 & & \textbf{7} \\
                6 & 5 & 8
            \end{tabular}
            &
            \begin{tabular}{l l l}
                1 & & \\
                2 & & 5 \\
                4 & 3 & \textbf{6}
            \end{tabular}
        \end{tabular}\\ \\

        \begin{tabular}{l p{56pt} p{72pt} l p{56pt} l}
            \begin{minipage}[c]{48pt}
                Add the \textbf{2} to middle pile:
            \end{minipage}
            &
            \begin{tabular}{l l l}
                1 & & \\
                4 & \textbf{2} & 7 \\
                6 & 5 & 8
            \end{tabular}
            &
            \begin{tabular}{l l l}
               1 & & \\
               2 & 3 & 5 \\
               4 & \textbf{7} & 6
            \end{tabular}
            &
            \begin{minipage}[c]{48pt}
                Add the \textbf{3} to right pile:
            \end{minipage}
            &
            \begin{tabular}{l l l}
                1 & & \textbf{3} \\
                4 & 2 & 7 \\
                6 & 5 & 8
            \end{tabular}
            &
            \begin{tabular}{l l l}
               1 & & 5 \\
               2 & 3 & 6 \\
               4 & 7 & \textbf{8}
            \end{tabular}
        \end{tabular}\\ \\

    \end{exgr}

    The idea behind Algorithm~\ref{alg:ExtendedPSalgorithm} is that we
    are using a new pile configuration $S(\sigma)$ (called the
    ``recording piles'') to implicitly label the order in which the
    elements of the permutation $\sigma$ are added to the usual
    Patience Sorting pile configuration $R(\sigma)$ (which we will now
    by analogy to RSK also call the ``insertion piles'').  It is clear
    that this information then allows us to uniquely reconstruct
    $\sigma$ by reversing the order in which the cards were played.
    However, even though reversing the Extended Patience Sorting
    Algorithm is much easier than reversing the RSK Algorithm through
    recursive ``reverse row bumping,'' the trade-off is that the pairs
    of pile configurations that result from the Extended Patience
    Sorting Algorithm are not independent (see \cite{refBLFPSAC05} for
    more details), whereas the standard Young tableau pairs generated
    by RSK are completely independent (up to shape).

    \section{Northeast Shadow Diagrams and Viennot's\\ Geometric RSK}
    \label{sec:GeometricRSK}

    In this section we briefly develop Viennot's geometric form for
    RSK in order to motivate the geometric form for the Extended
    Patience Sorting that is introduced in
    Section~\ref{sec:GeometricPS} below.

        \subsection{The Northeast Shadow Diagram of a Permutation}
        \label{sec:GeometricRSK:NEshadowDiagrams}

        We begin with the following fundamental definition:

        \begin{defn}\label{defn:Shadow}
            Given a lattice point $(m, n) \in \mathbb{Z}^{2}$, we
            define the \emph{northeast shadow} of $(m, n)$ to be the
            quarter space $S_{NE}(m, n) = \{ (x, y) \in \mathbb{R}^{2}
            \ | \ x \geq m, \ y \geq n\}$.
        \end{defn}

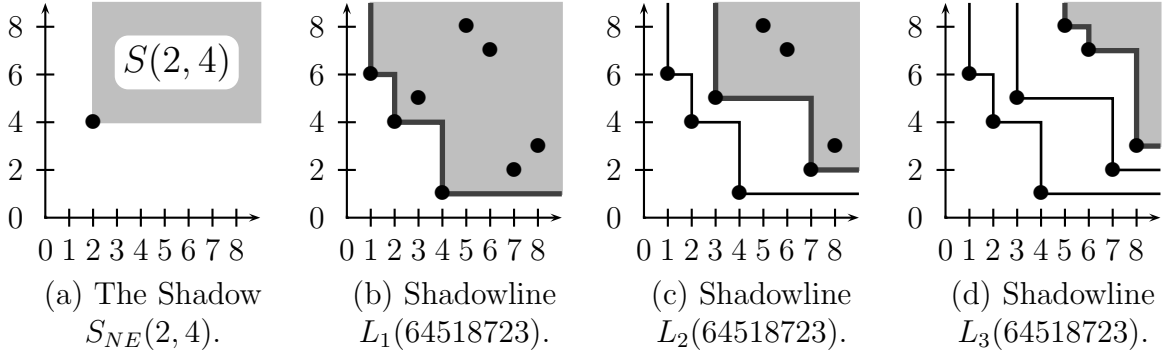
\begin{figure}[t]\label{fig:ShadowExample}
              \centering

              \begin{tabular}{cccc}

                  \begin{minipage}[c]{1.325in}
                      \centering

                          \psset{xunit=0.125in,yunit=0.125in}

                          \begin{pspicture}(0,0)(9,9)

                              \psaxes[Dy=2]{->}(9,9)

                              \pspolygon[fillcolor=lightgray, linecolor=lightgray,
                              fillstyle=solid](2,4)(9,4)(9,9)(2,9)

                              \rput(5.5,6.5){\psframebox*[framearc=.5]{{\large
                              $S(2,4)$}}}

                              \rput(2,4){{\large $\bullet$}}

                          \end{pspicture}\\[0.25in]

                        (a) The Shadow $S_{NE}(2,4)$.
                  \end{minipage}
                  &
                  \begin{minipage}[c]{1.5in}
                      \centering

                          \psset{xunit=0.125in,yunit=0.125in}

                          \begin{pspicture}(0,0)(9,9)

                              \psaxes[Dy=2]{->}(9,9)

                              \pspolygon[fillcolor=lightgray, linecolor=lightgray,
                              fillstyle=solid](1,9)(1,6)(2,6)(2,4)(4,4)(4,1)(9,1)(9,9)

                              \psline[linecolor=darkgray,linewidth=2pt](1,9)(1,6)(2,6)(2,4)(4,4)(4,1)(9,1)

                               \rput(1,6){{\large $\bullet$}}%
                               \rput(2,4){{\large $\bullet$}}%
                               \rput(3,5){{\large $\bullet$}}%
                               \rput(4,1){{\large $\bullet$}}%
                               \rput(5,8){{\large $\bullet$}}%
                               \rput(6,7){{\large $\bullet$}}%
                               \rput(7,2){{\large $\bullet$}}%
                               \rput(8,3){{\large $\bullet$}}%

                          \end{pspicture}\\[0.25in]

                      (b) Shadowline $L_{1}(64518723)$.
                  \end{minipage}

                  \begin{minipage}[c]{1.5in}
                      \centering

                          \psset{xunit=0.125in,yunit=0.125in}

                          \begin{pspicture}(0,0)(9,9)

                              \psaxes[Dy=2]{->}(9,9)

                              \pspolygon[fillcolor=lightgray, linecolor=lightgray,
                              fillstyle=solid](3,9)(3,5)(7,5)(7,2)(9,2)(9,9)

                              \psline[linecolor=darkgray,linewidth=2pt](3,9)(3,5)(7,5)(7,2)(9,2)

                              \psline[linecolor=black,linewidth=1pt](1,9)(1,6)(2,6)(2,4)(4,4)(4,1)(9,1)

                               \rput(1,6){{\large $\bullet$}}%
                               \rput(2,4){{\large $\bullet$}}%
                               \rput(3,5){{\large $\bullet$}}%
                               \rput(4,1){{\large $\bullet$}}%
                               \rput(5,8){{\large $\bullet$}}%
                               \rput(6,7){{\large $\bullet$}}%
                               \rput(7,2){{\large $\bullet$}}%
                               \rput(8,3){{\large $\bullet$}}%

                          \end{pspicture}\\[0.25in]

                      (c) Shadowline $L_{2}(64518723)$.
                  \end{minipage}
                  &
                  \begin{minipage}[c]{1.325in}
                      \centering

                          \psset{xunit=0.125in,yunit=0.125in}

                          \begin{pspicture}(0,0)(9,9)

                              \psaxes[Dy=2]{->}(9,9)

                              \pspolygon[fillcolor=lightgray, linecolor=lightgray,
                              fillstyle=solid](5,9)(5,8)(6,8)(6,7)(8,7)(8,3)(9,3)(9,9)

                              \psline[linecolor=darkgray,linewidth=2pt](5,9)(5,8)(6,8)(6,7)(8,7)(8,3)(9,3)

                              \psline[linecolor=black,linewidth=1pt](3,9)(3,5)(7,5)(7,2)(9,2)

                              \psline[linecolor=black,linewidth=1pt](1,9)(1,6)(2,6)(2,4)(4,4)(4,1)(9,1)

                               \rput(1,6){{\large $\bullet$}}%
                               \rput(2,4){{\large $\bullet$}}%
                               \rput(3,5){{\large $\bullet$}}%
                               \rput(4,1){{\large $\bullet$}}%
                               \rput(5,8){{\large $\bullet$}}%
                               \rput(6,7){{\large $\bullet$}}%
                               \rput(7,2){{\large $\bullet$}}%
                               \rput(8,3){{\large $\bullet$}}%

                          \end{pspicture}\\[0.25in]

                      (d) Shadowline $L_{3}(64518723)$.
                  \end{minipage}

              \end{tabular}\\

              \caption{Examples of Northeast Shadow and Shadowline
              Constructions}

          \end{figure}

        \noindent See Figure~\ref{fig:ShadowExample}(a) for an example
        of a point's northeast shadow.\\

        The most important use of these shadows is in building
        so-called northeast shadowlines:

        \begin{defn}\label{defn:Shadowline}
            Given lattice points $(m_{1}, n_{1}), (m_{2}, n_{2}),
            \ldots, (m_{k}, n_{k}) \in \mathbb{Z}^{2}$, we define
            their \emph{northeast shadowline} to be the boundary of
            the union of the northeast shadows $S_{NE}(m_{1}, n_{1}),
            S_{NE}(m_{2}, n_{2}), \ldots, S_{NE}(m_{k}, n_{k})$.
        \end{defn}

        In particular, we wish to associate to each permutation a
        certain collection of northeast shadowlines (as illustrated in
        Figure~\ref{fig:ShadowExample}(b)--(d)):

        \begin{defn}\label{defn:ShadowDiagram}
            Given a permutation $\sigma =
            \sigma_{1}\sigma_{2}\cdots\sigma_{n} \in
            \mathfrak{S}_{n}$, the \emph{northeast shadow diagram}
            $D_{NE}^{(0)}(\sigma)$ of $\sigma$ consists of the shadowlines
            $L_{1}(\sigma), L_{2}(\sigma), \ldots, L_{k}(\sigma)$
            formed as follows:\medskip

            \begin{itemize}
                \item $L_{1}(\sigma)$ is the northeast shadowline for
                the lattice points $\{(1, \sigma_{1}), (2,
                \sigma_{2}), \ldots, (n, \sigma_{n})\}$.\medskip

                \item While at least one of the points $(1,
                \sigma_{1}), (2, \sigma_{2}), \ldots, (n, \sigma_{n})$
                is not contained in the shadowlines $L_{1}(\sigma),
                L_{2}(\sigma), \ldots, L_{j}(\sigma)$, define
                $L_{j+1}(\sigma)$ to be the northeast shadowline for
                the points $$\{(i, \sigma_{i}) \ | \ (i, \sigma_{i})
                \notin \bigcup^{j}_{k=1} L_{k}(\sigma)\}.$$
            \end{itemize}

        \end{defn}

        In other words, we define the shadow diagram inductively by
        first taking $L_{1}(\sigma)$ to be the shadowline for the
        diagram $\{(1, \sigma_{1}), (2, \sigma_{2}), \ldots, (n,
        \sigma_{n})\}$ of the permutation.  Then we ignore the lattice
        points whose shadows were used in building $L_{1}(\sigma)$ and
        define $L_{2}(\sigma)$ to be the shadowline of the resulting
        subset of the permutation diagram.  We then build
        $L_{3}(\sigma)$ as the shadowline for the points not yet used
        in constructing either $L_{1}(\sigma)$ or $L_{2}(\sigma)$, and
        this process continues until all points in the permutation
        diagram are exhausted.

        We can characterize the points whose shadows define the
        shadowlines at each stage as follows: they are the smallest
        collection of unused points whose shadows collectively contain
        all other remaining unused points (and hence also contain the
        shadows of those points).  As a consequence of this shadow
        containment property, the shadowlines in a northeast shadow
        diagram will never cross.  However, as we will see in
        Section~\ref{sec:GeometricPS:SWshadowDiagrams} below, the dual
        construction to Definition~\ref{defn:ShadowDiagram} that is
        introduced will allow for crossing shadowlines, which are then
        classified and characterized in Section~\ref{sec:PSvsRSK}.
        The most fundamental cause for this distinction is the way
        that we will reverse the above shadow containment property for
        the points used in defining southwest shadowlines.

        \subsection{Viennot's Geometric RSK Algorithm}
        \label{sec:GeometricRSK:GeometricRSKalgorithm}

        \begin{figure}[t]\label{fig:GeometricRSKexample}
              \centering

              \begin{tabular}{cccc}

                  \begin{minipage}[c]{1.5in}
                      \centering

                          \psset{xunit=0.125in,yunit=0.125in}

                          \begin{pspicture}(0,0)(9,9)

                              \psaxes[Dy=2]{->}(9,9)

                              \psline[linecolor=black,linewidth=1pt](5,9)(5,8)(6,8)(6,7)(8,7)(8,3)(9,3)

                              \psline[linecolor=black,linewidth=1pt](3,9)(3,5)(7,5)(7,2)(9,2)

                              \psline[linecolor=black,linewidth=1pt](1,9)(1,6)(2,6)(2,4)(4,4)(4,1)(9,1)

                               \rput(1,6){{\large $\bullet$}}%
                               \rput(2,4){{\large $\bullet$}}%
                               \rput(3,5){{\large $\bullet$}}%
                               \rput(4,1){{\large $\bullet$}}%
                               \rput(5,8){{\large $\bullet$}}%
                               \rput(6,7){{\large $\bullet$}}%
                               \rput(7,2){{\large $\bullet$}}%
                               \rput(8,3){{\large $\bullet$}}%

                               \rput(2,6){{\large $\odot$}}%
                               \rput(4,4){{\large $\odot$}}%
                               \rput(7,5){{\large $\odot$}}%
                               \rput(6,8){{\large $\odot$}}%
                               \rput(8,7){{\large $\odot$}}%

                          \end{pspicture}\\[0.25in]

                        (a) Salient points for $D_{NE}^{(0)}(64518723)$.
                  \end{minipage}
                  &
                  \begin{minipage}[c]{1.5in}
                      \centering

                          \psset{xunit=0.125in,yunit=0.125in}

                          \begin{pspicture}(0,0)(9,9)

                              \psaxes[Dy=2]{->}(9,9)

                              \psline[linecolor=black,linewidth=1pt](9,7)(8,7)(8,9)

                              \psline[linecolor=black,linewidth=1pt](6,9)(6,8)(7,8)(7,5)(9,5)

                              \psline[linecolor=black,linewidth=1pt](2,9)(2,6)(4,6)(4,4)(9,4)

                               \rput(2,6){{\large $\bullet$}}%
                               \rput(4,4){{\large $\bullet$}}%
                               \rput(7,5){{\large $\bullet$}}%
                               \rput(6,8){{\large $\bullet$}}%
                               \rput(8,7){{\large $\bullet$}}%

                          \end{pspicture}\\[0.25in]

                      (b) Shadow Diagram $D_{NE}^{(1)}(64518723)$.
                  \end{minipage}

                  \begin{minipage}[c]{1.5in}
                      \centering

                          \psset{xunit=0.125in,yunit=0.125in}

                          \begin{pspicture}(0,0)(9,9)

                              \psaxes[Dy=2]{->}(9,9)

                              \psline[linecolor=black,linewidth=1pt](9,7)(8,7)(8,9)

                              \psline[linecolor=black,linewidth=1pt](6,9)(6,8)(7,8)(7,5)(9,5)

                              \psline[linecolor=black,linewidth=1pt](2,9)(2,6)(4,6)(4,4)(9,4)

                               \rput(2,6){{\large $\bullet$}}%
                               \rput(4,4){{\large $\bullet$}}%
                               \rput(7,5){{\large $\bullet$}}%
                               \rput(6,8){{\large $\bullet$}}%
                               \rput(8,7){{\large $\bullet$}}%

                               \rput(4,6){{\large $\odot$}}%
                               \rput(7,8){{\large $\odot$}}%

                          \end{pspicture}\\[0.25in]

                        (c) Salient points for $D_{NE}^{(1)}(64518723)$.
                  \end{minipage}
                  &
                  \begin{minipage}[c]{1.5in}
                      \centering

                          \psset{xunit=0.125in,yunit=0.125in}

                          \begin{pspicture}(0,0)(9,9)

                              \psaxes[Dy=2]{->}(9,9)

                              \psline[linecolor=black,linewidth=1pt](7,9)(7,8)(9,8)

                              \psline[linecolor=black,linewidth=1pt](4,9)(4,6)(9,6)

                               \rput(4,6){{\large $\bullet$}}%
                               \rput(7,8){{\large $\bullet$}}%

                          \end{pspicture}\\[0.25in]

                      (d) Shadow Diagram $D_{NE}^{(2)}(64518723)$.
                  \end{minipage}

              \end{tabular}\\

              \caption{The northeast shadow diagrams for the permutation $64518723 \in \mathfrak{S}_{8}$.}

          \end{figure}
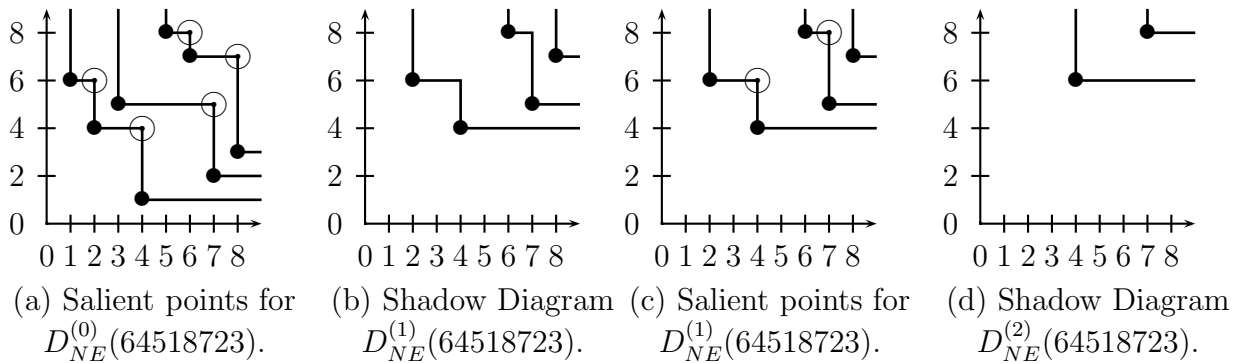

        As simple as northeast shadowlines were to define in the
        previous section, a great deal of information can still be
        gotten from them.  One of the most basic properties of the
        northeast shadow diagram $D_{NE}^{(0)}(\sigma)$ for a
        permutation $\sigma \in \mathfrak{S}_{n}$ is that it encodes
        the top row of the RSK insertion tableau $P(\sigma)$ (resp.
        recording tableau $Q(\sigma)$) as the smallest ordinates
        (resp.  smallest abscissae) of all points belonging to the
        constituent shadowlines $L_{1}(\sigma), L_{2}(\sigma), \ldots,
        L_{k}(\sigma)$.  One proves this by comparing the use of
        Schensted Insertion on the top row of the insertion tableau
        with the intersection of vertical lines having the form $x =
        a$.  In particular, as $a$ increases from $0$ to $n$, the line
        $x = a$ intersects the lattice points in the permutation
        diagram in the order that they are inserted into the top row,
        and so shadowlines connect elements of $\sigma$ to those
        smaller elements that will eventually bump them.  (See Sagan
        \cite{refSagan2000} for more details.)

        Remarkably, one can then use the northeast corners (called the
        \emph{salient points}) of $D_{NE}^{(0)}(\sigma)$ to form a new
        shadow diagram $D_{NE}^{(1)}(\sigma)$ that similarly gives the
        second rows of $P(\sigma)$ and $Q(\sigma)$.  Then,
        inductively, the salient points of $D_{NE}^{(1)}(\sigma)$ can
        be used to give the third rows of $P(\sigma)$ and $Q(\sigma)$,
        and so on.  As such, one can view this recursive formation of
        shadow diagrams as a geometric form for the RSK
        correspondence.  We illustrate this process in
        Figure~\ref{fig:GeometricRSKexample} for the following
        permutation from Example~\ref{eg:ExtendedPSexample}:

        \begin{displaymath}
            \sigma =  64518723 \stackrel{RSK}{\longleftrightarrow}
            \left(~
                \young(123,457,68)
                \raisebox{-0.5cm}{,}\
                \young(135,268,47)
            ~\right)
        \end{displaymath}

    \section{Southwest Shadow Diagrams and Geometric\\ Patience Sorting}
    \label{sec:GeometricPS}

    In this section we introduce a very natural dual algorithm to
    Viennot's geometric form for RSK as given in
    Section~\ref{sec:GeometricRSK:GeometricRSKalgorithm} above.

        \subsection{The Southwest Shadow Diagram of a Permutation}
        \label{sec:GeometricPS:SWshadowDiagrams}

        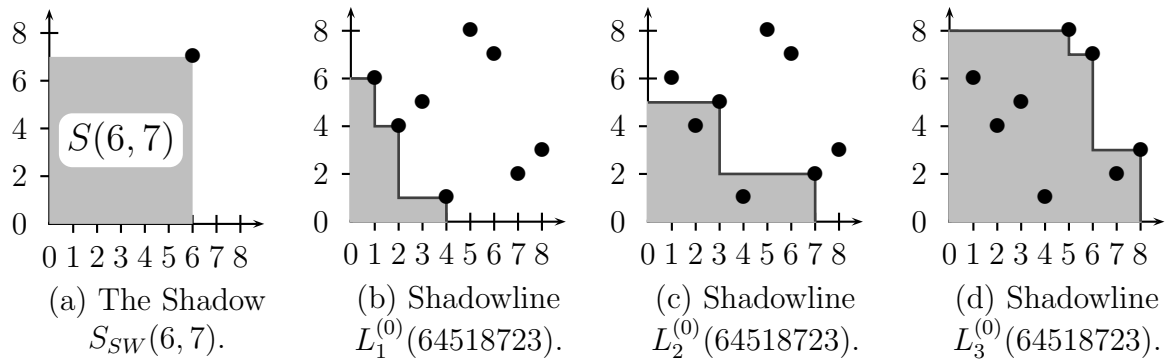
\begin{figure}[t]\label{fig:PSshadowExample}
              \centering

              \begin{tabular}{cccc}

                  \begin{minipage}[c]{1.325in}
                      \centering

                          \psset{xunit=0.125in,yunit=0.125in}

                          \begin{pspicture}(0,0)(9,9)

                              \psaxes[Dy=2]{->}(9,9)

                              \psframe[fillcolor=lightgray,fillstyle=solid,linecolor=lightgray](6,7)(0,0)

                              \rput(3,3.5){\psframebox*[framearc=.5]{{\large
                              $S(6,7)$}}}

                              \rput(6,7){{\large $\bullet$}}

                          \end{pspicture}\\[0.25in]

                        (a) The Shadow $S_{SW}(6,7)$.
                  \end{minipage}
                  &
                  \begin{minipage}[c]{1.5in}
                      \centering

                          \psset{xunit=0.125in,yunit=0.125in}

                        \begin{pspicture}(0,0)(9,9)

                            \psaxes[Dy=2]{->}(9,9)

                            \pspolygon[fillcolor=lightgray, linecolor=lightgray,
                            fillstyle=solid](0,6)(1,6)(1,4)(2,4)(2,1)(4,1)(4,0)(0,0)

                            \psline[linecolor=darkgray,
                            linewidth=1pt](0,6)(1,6)(1,4)(2,4)(2,1)(4,1)(4,0)

                             \rput(1,6){{\large $\bullet$}}%
                             \rput(2,4){{\large $\bullet$}}%
                             \rput(3,5){{\large $\bullet$}}%
                             \rput(4,1){{\large $\bullet$}}%
                             \rput(5,8){{\large $\bullet$}}%
                             \rput(6,7){{\large $\bullet$}}%
                             \rput(7,2){{\large $\bullet$}}%
                             \rput(8,3){{\large $\bullet$}}%

                        \end{pspicture}\\[0.25in]

                      (b) Shadowline $L_{1}^{(0)}(64518723)$.
                  \end{minipage}

                  \begin{minipage}[c]{1.5in}
                      \centering

                        \psset{xunit=0.125in,yunit=0.125in}

                        \begin{pspicture}(0,0)(9,9)

                            \psaxes[Dy=2]{->}(9,9)

                            \pspolygon[fillcolor=lightgray, linecolor=lightgray,
                            fillstyle=solid](0,5)(3,5)(3,2)(7,2)(7,0)(0,0)

                            \psline[linecolor=darkgray,
                            linewidth=1pt](0,5)(3,5)(3,2)(7,2)(7,0)

                             \rput(1,6){{\large $\bullet$}}%
                             \rput(2,4){{\large $\bullet$}}%
                             \rput(3,5){{\large $\bullet$}}%
                             \rput(4,1){{\large $\bullet$}}%
                             \rput(5,8){{\large $\bullet$}}%
                             \rput(6,7){{\large $\bullet$}}%
                             \rput(7,2){{\large $\bullet$}}%
                             \rput(8,3){{\large $\bullet$}}%

                        \end{pspicture}\\[0.25in]

                      (c) Shadowline $L_{2}^{(0)}(64518723)$.
                  \end{minipage}
                  &
                  \begin{minipage}[c]{1.325in}
                      \centering

                        \psset{xunit=0.125in,yunit=0.125in}

                        \begin{pspicture}(0,0)(9,9)

                            \psaxes[Dy=2]{->}(9,9)

                            \pspolygon[fillcolor=lightgray, linecolor=lightgray,
                            fillstyle=solid](0,8)(5,8)(5,7)(6,7)(6,3)(8,3)(8,0)(0,0)

                            \psline[linecolor=darkgray,
                            linewidth=1pt](0,8)(5,8)(5,7)(6,7)(6,3)(8,3)(8,0)

                             \rput(1,6){{\large $\bullet$}}%
                             \rput(2,4){{\large $\bullet$}}%
                             \rput(3,5){{\large $\bullet$}}%
                             \rput(4,1){{\large $\bullet$}}%
                             \rput(5,8){{\large $\bullet$}}%
                             \rput(6,7){{\large $\bullet$}}%
                             \rput(7,2){{\large $\bullet$}}%
                             \rput(8,3){{\large $\bullet$}}%

                        \end{pspicture}\\[0.25in]

                      (d) Shadowline $L_{3}^{(0)}(64518723)$.
                  \end{minipage}

              \end{tabular}\\

              \caption{Examples of Southwest Shadow and Shadowline
              Constructions}

          \end{figure}

        As in Section~\ref{sec:GeometricRSK:NEshadowDiagrams}, we
        begin with the following fundamental definition:

        \begin{defn}\label{defn:PSshadow}
            Given a lattice point $(m, n) \in \mathbb{Z}^{2}$, we
            define the \emph{southwest shadow} of $(m, n)$ to be the
            quarter space $S_{SW}(m, n) = \{ (x, y) \in \mathbb{R}^{2} \ |
            \ x \leq m, \ y \leq n\}$.
        \end{defn}

        \noindent See Figure~\ref{fig:PSshadowExample}(a) for an
        example of a point's southwest shadow.\\

        As with their northeast counterparts, the most important use
        of these shadows is in building southwest shadowlines:

        \begin{defn}\label{defn:PSshadowline}
            Given lattice points $(m_{1}, n_{1}), (m_{2}, n_{2}),
            \ldots, (m_{k}, n_{k}) \in \mathbb{Z}^{2}$, we define
            their \emph{southwest shadowline} to be the boundary of
            the union of the shadows $S_{SW}(m_{1}, n_{1})$,
            $S_{SW}(m_{2}, n_{2})$, $\ldots$, $S_{SW}(m_{k}, n_{k})$.
        \end{defn}

        In particular, we wish to associate to each permutation a
        certain collection of southwest shadowlines.  However, unlike
        the northeast case, these shadowlines sometimes cross (as
        illustrated in Figures~\ref{fig:PSshadowExample}(b)--(d) and
        Figure~\ref{fig:GeometricPS_example}(a)).

        \begin{defn}\label{defn:PSshadowDiagram}
            Given a permutation $\sigma =
            \sigma_{1}\sigma_{2}\cdots\sigma_{n} \in
            \mathfrak{S}_{n}$, the \emph{southwest shadow diagram}
            $D_{SW}^{(0)}(\sigma)$ of $\sigma$ consists of the
            southwest shadowlines $L_{1}^{(0)}(\sigma),
            L_{2}^{(0)}(\sigma), \ldots, L_{k}^{(0)}(\sigma)$ formed
            as follows:\medskip

            \begin{itemize}
                \item $L_{1}^{(0)}(\sigma)$ is the shadowline for
                those lattice points $(x, y) \in \{(1, \sigma_{1}),
                (2, \sigma_{2}), \ldots, (n, \sigma_{n})\}$ such that
                $S_{SW}(x, y)$ does not contain any other lattice
                points.\medskip

                \item While at least one of the points $(1,
                \sigma_{1}), (2, \sigma_{2}), \ldots, (n, \sigma_{n})$
                is not contained in the shadowlines
                $L_{1}^{(0)}(\sigma), L_{2}^{(0)}(\sigma), \ldots,
                L_{j}^{(0)}(\sigma)$, define $L_{j+1}^{(0)}(\sigma)$
                to be the shadowline for the points
                \[
                (x, y) \in \{(i, \sigma_{i}) \ | \ (i, \sigma_{i})
                \notin \bigcup^{j}_{k=1} L_{k}^{(0)}(\sigma)\}
                \]
                such that $S_{SW}(x, y)$ does not contain any other
                lattice points in the same set.
            \end{itemize}

        \end{defn}

        In other words, we again define a shadow diagram by
        recursively eliminating certain points in the permutation
        diagram until every point has been used to define a
        shadowline.  However, we are here reversing both the direction
        of the shadows and the shadow containment property from the
        northeast case.  It is in this sense that the geometric form
        for the Extended Patience Sorting Algorithm given in the next
        section can be viewed as ``dual'' to Viennot's geometric form
        for RSK.

        \subsection{The Geometric Patience Sorting Algorithm}
        \label{sec:GeometricPS:GeometricPSalgorithm}

        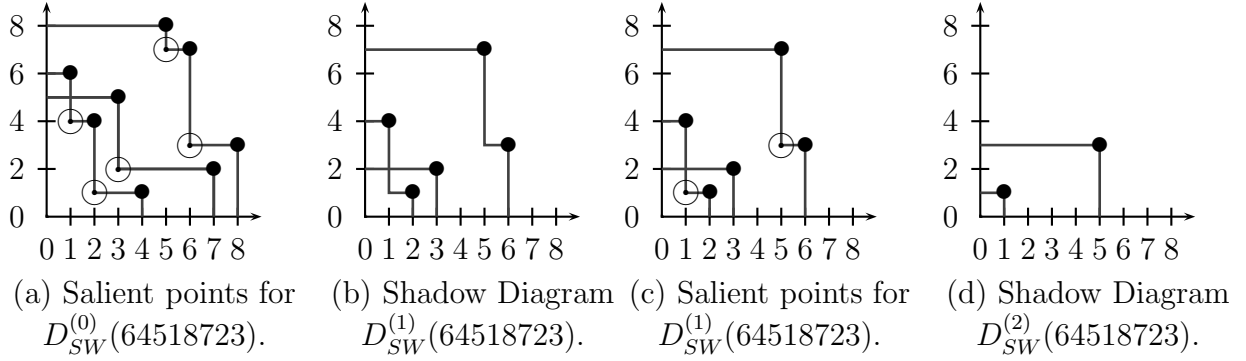
\begin{figure}[t]\label{fig:GeometricPS_example}
              \centering

              \begin{tabular}{cccc}

                  \begin{minipage}[c]{1.5in}
                      \centering

                        \psset{xunit=0.125in,yunit=0.125in}

                        \begin{pspicture}(0,0)(9,9)

                            \psline[linecolor=darkgray,
                            linewidth=1pt](0,5)(3,5)(3,2)(7,2)(7,0)

                            \psaxes[Dy=2]{->}(9,9)

                            \psline[linecolor=darkgray,
                            linewidth=1pt](0,6)(1,6)(1,4)(2,4)(2,1)(4,1)(4,0)

                            \psline[linecolor=darkgray,
                            linewidth=1pt](0,5)(3,5)(3,2)(7,2)(7,0)

                            \psline[linecolor=darkgray,
                            linewidth=1pt](0,8)(5,8)(5,7)(6,7)(6,3)(8,3)(8,0)

                            \rput(1,6){{\large $\bullet$}}%
                            \rput(2,4){{\large $\bullet$}}%
                            \rput(3,5){{\large $\bullet$}}%
                            \rput(4,1){{\large $\bullet$}}%
                            \rput(5,8){{\large $\bullet$}}%
                            \rput(6,7){{\large $\bullet$}}%
                            \rput(7,2){{\large $\bullet$}}%
                            \rput(8,3){{\large $\bullet$}}%

                            \rput(1,4){{\large $\odot$}}%
                            \rput(2,1){{\large $\odot$}}%
                            \rput(3,2){{\large $\odot$}}%
                            \rput(5,7){{\large $\odot$}}%
                            \rput(6,3){{\large $\odot$}}%

                        \end{pspicture}\\[0.25in]

                        (a) Salient points for $D_{SW}^{(0)}(64518723)$.
                  \end{minipage}
                  &
                  \begin{minipage}[c]{1.5in}
                      \centering

                        \psset{xunit=0.125in,yunit=0.125in}

                        \begin{pspicture}(0,0)(9,9)

                            \psaxes[Dy=2]{->}(9,9)

                            \psline[linecolor=darkgray,
                            linewidth=1pt](0,4)(1,4)(1,1)(2,1)(2,0)

                            \psline[linecolor=darkgray,
                            linewidth=1pt](0,2)(3,2)(3,0)

                            \psline[linecolor=darkgray,
                            linewidth=1pt](0,7)(5,7)(5,3)(6,3)(6,0)

                             \rput(1,4){{\large $\bullet$}}%
                             \rput(2,1){{\large $\bullet$}}%
                             \rput(3,2){{\large $\bullet$}}%
                             \rput(5,7){{\large $\bullet$}}%
                             \rput(6,3){{\large $\bullet$}}%

                        \end{pspicture}\\[0.25in]

                      (b) Shadow Diagram $D_{SW}^{(1)}(64518723)$.
                  \end{minipage}

                  \begin{minipage}[c]{1.5in}
                      \centering

                        \psset{xunit=0.125in,yunit=0.125in}

                        \begin{pspicture}(0,0)(9,9)

                            \psaxes[Dy=2]{->}(9,9)

                            \psline[linecolor=darkgray,
                            linewidth=1pt](0,4)(1,4)(1,1)(2,1)(2,0)

                            \psline[linecolor=darkgray,
                            linewidth=1pt](0,2)(3,2)(3,0)

                            \psline[linecolor=darkgray,
                            linewidth=1pt](0,7)(5,7)(5,3)(6,3)(6,0)

                             \rput(1,4){{\large $\bullet$}}%
                             \rput(2,1){{\large $\bullet$}}%
                             \rput(3,2){{\large $\bullet$}}%
                             \rput(5,7){{\large $\bullet$}}%
                             \rput(6,3){{\large $\bullet$}}%

                             \rput(1,1){{\large $\odot$}}%
                             \rput(5,3){{\large $\odot$}}%

                        \end{pspicture}\\[0.25in]

                        (c) Salient points for $D_{SW}^{(1)}(64518723)$.
                  \end{minipage}
                  &
                  \begin{minipage}[c]{1.5in}
                      \centering

                        \psset{xunit=0.125in,yunit=0.125in}

                        \begin{pspicture}(0,0)(9,9)

                            \psaxes[Dy=2]{->}(9,9)

                            \psline[linecolor=darkgray,
                            linewidth=1pt](0,1)(1,1)(1,0)

                            \psline[linecolor=darkgray,
                            linewidth=1pt](0,3)(5,3)(5,0)

                             \rput(1,1){{\large $\bullet$}}%
                             \rput(5,3){{\large $\bullet$}}%

                        \end{pspicture}\\[0.25in]

                      (d) Shadow Diagram $D_{SW}^{(2)}(64518723)$.
                  \end{minipage}

              \end{tabular}\\

              \caption{The southwest shadow diagrams for the permutation $64518723 \in \mathfrak{S}_{8}$.}

          \end{figure}

        As in Section~\ref{sec:GeometricRSK:GeometricRSKalgorithm},
        one can produce a sequence $D_{SW}(\sigma) =
        (D_{SW}^{(0)}(\sigma), D_{SW}^{(1)}(\sigma),
        D_{SW}^{(2)}(\sigma), \ldots)$ of shadow diagrams for a given
        permutation $\sigma \in \mathfrak{S}_{n}$ by recursively
        applying Definition~\ref{defn:PSshadowDiagram} to salient
        points, with the restriction that new shadowlines can only
        connect points that were on the same shadowline in the
        previous iteration.  (The reason for this important
        distinction from Geometric RSK is discussed further in
        Section~\ref{sec:PSvsRSK:TypesOfCrossings} below.)  The
        salient points in this case are then naturally defined to be
        the southwest corner points of a given set of shadowlines.
        See Figure~\ref{fig:GeometricPS_example} for an example of how
        this works for the permutation from
        Example~\ref{eg:ExtendedPSexample}.

        \begin{defn}

            We call $D_{SW}^{(k)}(\sigma)$ the $k^{\rm th}$ \emph{iterate}
            of the \emph{exhaustive shadow diagram} $D_{SW}(\sigma)$ for
            $\sigma \in \mathfrak{S}_{n}$.

        \end{defn}

        Moreover, the resulting sequence of shadow diagrams can then
        be used to reconstruct the pair of pile configurations given
        by the Extended Patience Sorting Algorithm
        (Algorithm~\ref{alg:ExtendedPSalgorithm}).  To accomplish
        this, index the cards in a pile configuration using the French
        convention for tableaux so that the row index increases from
        bottom to top and the column index from left to right.  (In
        other words, we are labelling boxes as we would lattice points
        in the first quadrant of $\mathbb{R}^2$).  Then, for a given
        permutation $\sigma \in \mathfrak{S}_{n}$, the elements of the
        $i$th row of the insertion piles $R(\sigma)$ (resp.  recording
        piles $S(\sigma)$) are given by the largest ordinates (resp.
        abscissae) of the shadowlines that compose $D_{SW}^{(i)}$.

        The main difference between this process and Viennot's
        Geometric RSK is that care must be taken to assemble each row
        in its proper order.  Unlike the entries of a Young tableau,
        the elements in the rows of a pile configuration do not
        necessarily increase from left to right, and they do not have
        to be contiguous.  As such, the components of each row should
        be recorded in the order that the shadowlines are formed.  The
        rows can then uniquely be assembled into a legal pile
        configuration since the elements in the columns of a pile
        configuration must both decrease (when read from bottom to
        top) and appear in the leftmost pile possible.\\

        The proof of this is along the same lines as that of Viennot's
        Geometric RSK in that the shadowlines produced by
        Definition~\ref{defn:PSshadowDiagram} are a visual record for
        how cards are played atop each other under
        Algorithm~\ref{alg:ExtendedPSalgorithm}.  In particular, it
        should be clear that, given a permutation $\sigma \in
        \mathfrak{S}_{n}$, the shadowlines in both of the shadow
        diagrams $D_{SW}^{(0)}(\sigma)$ and $D_{NE}^{(0)}(\sigma)$ are
        defined by the same lattice points from the permutation
        diagram for $\sigma$.  In \cite{refBLFPSAC05} the points along
        a given northeast shadowline are shown to correspond exactly
        to the elements in some column of $R(\sigma)$ (as both
        correspond to one of the left-to-right minima subsequences of
        $\sigma$).  Thus, by reading the lattice points in the
        permutation diagram in increasing order of their abscissae,
        one can uniquely reconstruct both the piles in $R(\sigma)$ and
        the exact order in which cards are added to these piles (which
        implicitly yields $S(\sigma)$).  In this sense, both
        $D_{SW}^{(0)}(\sigma)$ and $D_{NE}^{(0)}(\sigma)$ encode the
        bottom rows of $R(\sigma)$ and $S(\sigma)$ as given by
        Algorithm~\ref{alg:ExtendedPSalgorithm}.

        It is then easy to see by induction that the salient points of
        $D_{SW}^{(k-1)}(\sigma)$ yield the $k^{\rm th}$ rows of
        $R(\sigma)$ and $S(\sigma)$, and so this gives the following

        \begin{thm}
            The process described above for creating a pair of pile
            configurations $(R'(\sigma), S'(\sigma))$ from the
            Geometric Patience Sorting construction yields the same
            pair of pile configurations $(R(\sigma), S(\sigma))$ as
            the Extended Patience Sorting Algorithm
            (Algorithm~\ref{alg:ExtendedPSalgorithm}).
        \end{thm}

        Having given the above Geometric form for
        Algorithm~\ref{alg:ExtendedPSalgorithm}, it is worth pointing
        out that---as with RSK---there are various natural
        generalizations of Extended Patience Sorting for more general
        combinatorial objects including words and lexicographic
        arrays.  (See \cite{refFulton1997} for a description of such
        extensions of RSK.) Moreover, many of these generalizations
        can still similarly be realized as non-recursive analogs for
        the forms of RSK that can be applied to such objects.  In
        particular, the authors in \cite{refBLinPrep2005} explore
        several such generalizations and develop geometric forms for
        them much like the one given in this section.

        In the case of words, Aldous and Diaconis \cite{refAD1999}
        have given two different generalizations for Patience Sorting
        based upon whether cards with equal value are played on top of
        each other or not.  These are called the ``ties allowed'' and
        ``ties forbidden'' cases, respectively, and the usual RSK and
        dual RSK algorithms can be modeled in order to develop
        bijective versions of them.  The geometric forms for the
        resulting algorithms as given in \cite{refBLinPrep2005} can
        then be compared to Fulton's ``Matrix-Ball'' Geometric RSK
        algorithm (defined in \cite{refFulton1997}) just as we compare
        the Geometric Patience Sorting given in this section to
        Viennot's Geometric RSK in Section~\ref{sec:PSvsRSK}.

    \section{Geometric Patience Sorting and Intersecting\\ Lattice Paths}
    \label{sec:PSvsRSK}

        Extended Patience Sorting
        (Algorithm~\ref{alg:ExtendedPSalgorithm}) can be viewed as a
        ``non-bumping'' version of the RSK algorithm for permutations
        in that cards are permanently placed into piles and are
        covered by other cards rather being displaced by them.  It is
        in this sense that one of the main differences between their
        geometric algorithms lies in how and in what order (when read
        from left to right) the salient points of their respective
        shadow diagrams are determined.  In particular, as playing a
        card atop a pre-existing pile under Patience Sorting is
        essentially like non-recursive Schensted Insertion, certain
        particularly egregious ``multiple bumps'' that occur under the
        Schensted Insertion Algorithm prove to be too complicated to
        be properly modeled by the ``static insertions'' of Patience
        Sorting.

        At the same time, it is also easy to see that for a given
        $\sigma \in \mathfrak{S}_{n}$, the cards atop the piles in the
        pile configurations $R(\sigma)$ and $S(\sigma)$ (as given by
        Algorithm~\ref{alg:ExtendedPSalgorithm}) are exactly the cards
        in the top rows of the RSK insertion tableau $P(\sigma)$ and
        recording tableau $Q(\sigma)$, respectively.  Thus, this
        raises the question of when the remaining rows of $P(\sigma)$
        and $Q(\sigma)$ can likewise be recovered from $R(\sigma)$ and
        $S(\sigma)$.  While this appears to be directly related to the
        order in which salient points are read (as illustrated in
        Example~\ref{eg:2431RSKvsPSexample} below), one would
        ultimately hope to characterize the answer in terms of
        generalized pattern avoidance similar to the description of
        reverse patience words for pile configurations (as given in
        \cite{refBLFPSAC05}).

        \begin{exgr}\label{eg:2431RSKvsPSexample}
            Consider the northeast and southwest shadow diagrams for
            $\sigma = 2431$:\\ 

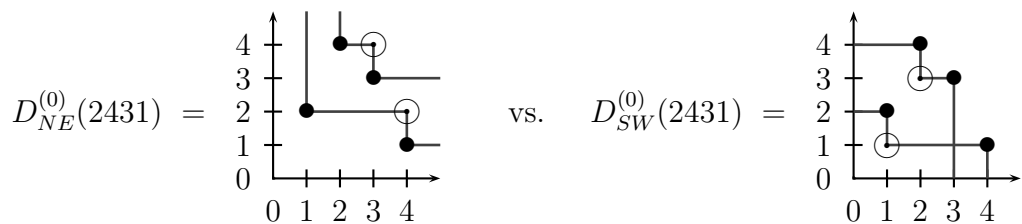
            \begin{center}

                \begin{tabular}{ccccc}

                    \raisebox{-.25cm}{$D_{NE}^{(0)}(2431) \ = $}
                    &
                    \begin{minipage}[c]{1.25in}
                        \centering

                          \psset{xunit=0.175in,yunit=0.175in}

                          \begin{pspicture}(0,0)(5,5)

                              \psaxes{->}(5,5)

                              \psline[linecolor=darkgray,
                              linewidth=1pt](1,5)(1,2)(4,2)(4,1)(5,1)

                              \psline[linecolor=darkgray,
                              linewidth=1pt](2,5)(2,4)(3,4)(3,3)(5,3)

                              \rput(1,2){{\large $\bullet$}}%
                              \rput(2,4){{\large $\bullet$}}%
                              \rput(3,3){{\large $\bullet$}}%
                              \rput(4,1){{\large $\bullet$}}%

                              \rput(4,2){{\large $\odot$}}%
                              \rput(3,4){{\large $\odot$}}%

                          \end{pspicture}
                    \end{minipage}
                    &
                    \raisebox{-.25cm}{vs. \quad $D_{SW}^{(0)}(2431) \ = $}
                    &
                    \begin{minipage}[c]{1.25in}
                        \centering

                          \psset{xunit=0.175in,yunit=0.175in}

                          \begin{pspicture}(0,0)(5,5)

                              \psaxes{->}(5,5)

                              \psline[linecolor=darkgray,
                              linewidth=1pt](0,2)(1,2)(1,1)(4,1)(4,0)

                              \psline[linecolor=darkgray,
                              linewidth=1pt](0,4)(2,4)(2,3)(3,3)(3,0)

                              \rput(1,2){{\large $\bullet$}}%
                              \rput(2,4){{\large $\bullet$}}%
                              \rput(3,3){{\large $\bullet$}}%
                              \rput(4,1){{\large $\bullet$}}%

                              \rput(1,1){{\large $\odot$}}%
                              \rput(2,3){{\large $\odot$}}%

                          \end{pspicture}
                    \end{minipage}\\[.65in]

                \end{tabular}
            \end{center}

            \noindent In particular, note that the order in which the
            salient points are formed (when read from left to right)
            is reversed.  Such reversals serve to illustrate one of
            the inherent philosophical differences between RSK and the
            Extended Patience Sorting Algorithm.

            As mentioned in
            Section~\ref{sec:GeometricPS:GeometricPSalgorithm} above,
            another fundamental difference between Geometric RSK and
            Geometric Patience Sorting is that the latter allows
            certain crossings to occur in the lattice paths formed
            during the same iteration of the algorithm.  We classify
            these crossings in
            Section~\ref{sec:PSvsRSK:TypesOfCrossings} and then
            characterize those permutations that yield entirely
            non-intersecting lattice paths in
            Section~\ref{sec:PSvsRSK:NoCrossings}.

        \end{exgr}

        \subsection{Types of Crossings in Geometric Patience Sorting}
        \label{sec:PSvsRSK:TypesOfCrossings}

            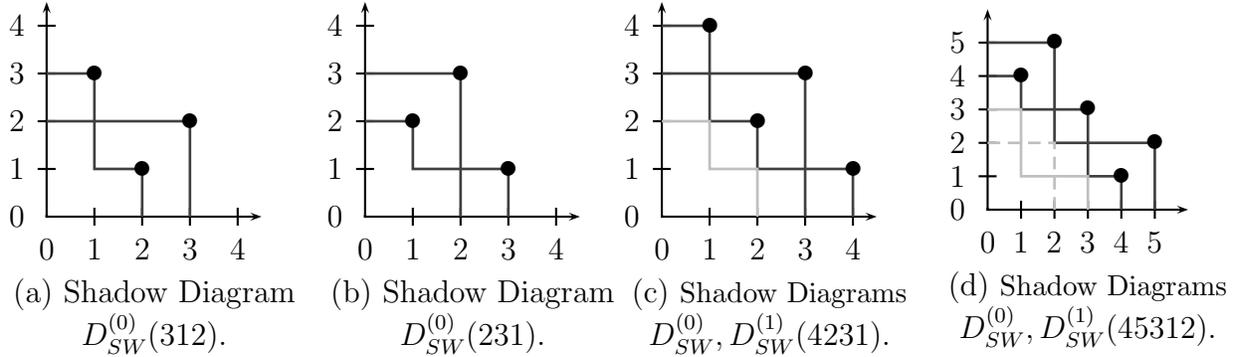
\begin{figure}[t]\label{fig:PSvsRSKexamples}
                  \centering

                  \begin{tabular}{cccc}

                      \begin{minipage}[c]{1.5in}
                          \centering

                            \psset{xunit=0.125in,yunit=0.125in}

                            \begin{pspicture}(0,0)(9,9)

                                \psaxes[dx=2,dy=2]{->}(9,9)

                                \psline[linecolor=darkgray,
                                linewidth=1pt](0,6)(2,6)(2,2)(4,2)(4,0)

                                \psline[linecolor=darkgray,
                                linewidth=1pt](0,4)(6,4)(6,0)

                                \rput(2,6){{\large $\bullet$}}%
                                \rput(4,2){{\large $\bullet$}}%
                                \rput(6,4){{\large $\bullet$}}%

                            \end{pspicture}\\[0.25in]

                            (a) Shadow Diagram $D_{SW}^{(0)}(312)$.
                      \end{minipage}
                      &
                      \begin{minipage}[c]{1.5in}
                          \centering

                            \psset{xunit=0.125in,yunit=0.125in}

                            \begin{pspicture}(0,0)(9,9)

                                \psaxes[dx=2,dy=2]{->}(9,9)

                                \psline[linecolor=darkgray,
                                linewidth=1pt](0,4)(2,4)(2,2)(6,2)(6,0)

                                \psline[linecolor=darkgray,
                                linewidth=1pt](0,6)(4,6)(4,0)

                                \rput(2,4){{\large $\bullet$}}%
                                \rput(4,6){{\large $\bullet$}}%
                                \rput(6,2){{\large $\bullet$}}%

                            \end{pspicture}\\[0.25in]

                            (b) Shadow Diagram $D_{SW}^{(0)}(231)$.
                      \end{minipage}

                      \begin{minipage}[c]{1.5in}
                          \centering

                            \psset{xunit=0.125in,yunit=0.125in}

                            \begin{pspicture}(0,0)(9,9)

                                \psaxes[dx=2,dy=2]{->}(9,9)

                                \psline[linecolor=darkgray,
                                linewidth=1pt](0,8)(2,8)(2,4)(4,4)(4,2)(8,2)(8,0)

                                \psline[linecolor=darkgray,
                                linewidth=1pt](0,6)(6,6)(6,0)

                                \psline[linecolor=lightgray, linestyle=solid,
                                linewidth=1pt](0,4)(2,4)(2,2)(4,2)(4,0)

                                \rput(2,8){{\large $\bullet$}}%
                                \rput(4,4){{\large $\bullet$}}%
                                \rput(6,6){{\large $\bullet$}}%
                                \rput(8,2){{\large $\bullet$}}%

                            \end{pspicture}\\[0.25in]

                            (c) {\small Shadow Diagrams} $D_{SW}^{(0)},
                            D_{SW}^{(1)}(4231)$.
                      \end{minipage}
                      &
                      \begin{minipage}[c]{1.5in}
                          \centering

                            \psset{xunit=0.175in,yunit=0.175in}

                            \begin{pspicture}(0,0)(6,6)

                                \psaxes{->}(6,6)

                                \psline[linecolor=darkgray,
                                linewidth=1pt](0,4)(1,4)(1,3)(3,3)(3,1)(4,1)(4,0)

                                \psline[linecolor=darkgray,
                                linewidth=1pt](0,5)(2,5)(2,2)(5,2)(5,0)

                                \psline[linecolor=lightgray, linestyle=solid,
                                linewidth=1pt](0,3)(1,3)(1,1)(3,1)(3,0)

                                \psline[linecolor=lightgray, linestyle=dashed,
                                linewidth=1pt](0,2)(2,2)(2,0)

                                \rput(1,4){{\large $\bullet$}}%
                                \rput(2,5){{\large $\bullet$}}%
                                \rput(3,3){{\large $\bullet$}}%
                                \rput(4,1){{\large $\bullet$}}%
                                \rput(5,2){{\large $\bullet$}}%

                            \end{pspicture}\\[0.25in]

                            (d) {\small Shadow Diagrams} $D_{SW}^{(0)},
                            D_{SW}^{(1)}(45312)$.
                      \end{minipage}

                  \end{tabular}\\

                  \caption{Shadow diagrams with different types of crossings.}

              \end{figure}

            Given $\sigma \in \mathfrak{S}_{n}$, we can classify the
            basic types of crossings in $D_{SW}^{(0)}(\sigma)$ as
            follows: First note that each southwest shadowline in
            $D_{SW}^{(0)}(\sigma)$ corresponds to a pair of decreasing
            sequences of the same length, namely a column from the
            insertion piles $R(\sigma)$ and its corresponding column
            from the recording piles $S(\sigma)$.  Then, given two
            different pairs of such columns in $R(\sigma)$ and
            $S(\sigma)$, the shadowline corresponding to the rightmost
            (resp.  leftmost) pair---under the convention that new
            columns are always added to the right of all other columns
            in Algorithm~\ref{alg:ExtendedPSalgorithm})---is called
            the \emph{upper} (resp.  \emph{lower}) shadowline. More
            formally:

            \begin{defn}
                Given two shadowlines, $L^{(m)}_{i}(\sigma),
                L^{(m)}_{j}(\sigma)\in D_{SW}^{(m)}(\sigma)$ with $i<j$, we
                call $L^{(m)}_{i}(\sigma)$ the \emph{lower} shadowline and
                $L^{(m)}_{j}(\sigma)$, the \emph{upper} shadowline.
                Moreover, if $L^{(m)}_{i}(\sigma)$ and
                $L^{(m)}_{j}(\sigma)$ intersect, then we call this a
                \emph{vertical crossing} (resp.  \emph{horizontal
                crossing}) if it involves a vertical (resp.  horizontal)
                segment of $L^{(m)}_{j}(\sigma)$.
            \end{defn}\medskip

               We illustrate these crossings in the following example.
               In particular, note that the only permutations $\sigma
               \in \mathfrak{S}_{3}$ of length three having
               intersections in their $0^{\textrm{th}}$ iterate shadow
               diagram $D_{SW}^{(0)}(\sigma)$ are $312, 231 \in
               \mathfrak{S}_{3}$.

            \begin{exgr} $\phantom{a}$

                \begin{enumerate}
                    \item The smallest permutation for which
                    $D_{SW}^{(0)}(\sigma)$ contains a horizontal
                    crossing is $\sigma = 312$ as illustrated in
                    Figure~\ref{fig:PSvsRSKexamples}(a).  The upper
                    shadowline involved in this crossing is the one
                    with only two segments.\smallskip

                    \item The smallest permutation for which
                    $D_{SW}^{(0)}(\sigma)$ contains a vertical
                    crossing is $\sigma = 231$ as illustrated in
                    Figure~\ref{fig:PSvsRSKexamples}(b).  As in part
                    (1), the upper shadowline involved in this
                    crossing is again the one with only two
                    segments.\smallskip

                    \item Consider $\sigma = 4231$.  From
                    Figure~\ref{fig:PSvsRSKexamples}(c),
                    $D_{SW}^{(0)}(\sigma)$ contains exactly two
                    southwest shadowlines, and these shadowlines form
                    a horizontal crossing followed by a vertical
                    crossing.  We call a configuration like this a
                    ``polygonal crossing.''  Note in particular that
                    $D_{SW}^{(1)}(\sigma)$ (trivially) has no
                    crossings.\smallskip

                    \item Consider $\sigma = 45312$.  From
                    Figure~\ref{fig:PSvsRSKexamples}(d),
                    $D_{SW}^{(0)}(\sigma)$ not only has a ``polygonal
                    crossing'' (this time as two shadowlines have a
                    vertical crossing followed by a horizontal one)
                    but $D_{SW}^{(1)}(\sigma)$ does as well.\smallskip

                \end{enumerate}

            \end{exgr}

            \noindent Polygonal crossings are what make it
            necessary to read only the salient points along the
            same shadowline in the order in which shadowlines are
            formed (as opposed to constructing the subsequent
            shadowlines using the entire partial permutation of
            salient points as in Viennot's Geometric RSK).

           \begin{exgr}

               Consider the shadow diagram of $\sigma=45312$ as
               illustrated in Figure~\ref{fig:PSvsRSKexamples}(d).
               The $0^{\textrm{th}}$ iterate shadow diagram
               $D_{SW}^{(0)}$ contain a polygonal crossing, and so
               the $1^{\textrm{st}}$ iterate shadow diagram
               $D_{SW}^{(1)}$ needs to be formed as indicated in
               order to properly describe the pile configurations
               $R(\sigma)$ and $S(\sigma)$ since

                   \begin{displaymath}
                       \sigma =  45312 \stackrel{XPS}{\longleftrightarrow}
                       \left(~
                           \begin{minipage}[c]{32pt}
                               $\begin{array}{cc}
                                   1 &   \\
                                   3 & 2  \\
                                   4 & 5
                               \end{array}$
                           \end{minipage}
                           \raisebox{-0.5cm}{,}\
                           \begin{minipage}[c]{32pt}
                               $\begin{array}{cc}
                                   1 &   \\
                                   3 & 2  \\
                                   4 & 5
                               \end{array}$
                           \end{minipage}
                       ~\right)
                   \end{displaymath}\medskip

               \noindent under the Extended Patience Sorting
               Algorithm.

           \end{exgr}

        \subsection{Non-intersecting shadow diagrams}
        \label{sec:PSvsRSK:NoCrossings}

            Unlike the rows of Young tableaux, the values in the rows
            of a pile configuration do not necessarily increase when
            read from left to right.  In fact, the descents in the
            rows of pile configurations are very closely related to
            the crossings given by Geometric Patience Sorting.
            
            As noted in
            Section~\ref{sec:GeometricPS:GeometricPSalgorithm} above,
            Geometric Patience Sorting is ostensibly simpler than
            Geometric RSK in that one can essentially recover both the
            insertion piles $R(\sigma)$ and the recording piles
            $S(\sigma)$ from the $0^{\textrm{th}}$ iterate shadow
            diagram $D_{SW}^{(0)}$.  The fundamental use, then, of the
            iterates $D_{SW}^{(i+1)}, D_{SW}^{(i+2)}, \ldots$ is in
            understanding the intersections in the $i^{\textrm{th}}$
            iterate shadow diagram $D_{SW}^{(i)}$.  In particular,
            each shadowline $L^{(m)}_i(\sigma)\in
            D_{SW}^{(m)}(\sigma)$ corresponds to the pair of segments
            of the $i^{\rm th}$ columns of $R(\sigma)$ and $S(\sigma)$
            that are above the $m^{\rm th}$ row (or are the $i^{\rm
            th}$ columns if $m=0$), where rows are numbered from
            bottom to top.

            \begin{thm}\label{thm:NoncrossingPilesCondition}
                Each iterate $D_{SW}^{(m)}(\sigma)$ ($m\ge 0$) of $\sigma \in
                \mathfrak{S}_{n}$ is free from crossings if and only if
                every row in both $R(\sigma)$ and $S(\sigma)$ is monotone
                increasing from left to right.
            \end{thm}

            \begin{proof}

                Since each $L^{(m)}_i=L^{(m)}_i(\sigma)$ depends only on
                the $i^{\rm th}$ columns of $R=R(\sigma)$ and
                $S=S(\sigma)$ above row $m$, we may assume without loss of
                generality that $R$ and $S$ have the same shape with
                exactly two columns.

                Let $m+1$ be the highest row where a descent occurs in
                either $R$ or $S$.  If this descent occurs in $R$, then
                $L^{(m)}_2$ is the upper shadowline in a horizontal
                crossing since $L^{(m)}_2$ has $y$-intercept below that of
                $L^{(m)}_1$, which is the lower shadowline in this
                crossing (as in $312$).  If this descent occurs in $S$,
                then $L^{(m)}_2$ is the upper shadowline in a vertical
                crossing since $L^{(m)}_2$ has $x$-intercept to the left
                of $L^{(m)}_1$, which is the lower shadowline in
                this crossing (as in $231$).  Note that both descents may
                occur simultaneously (as in $4231$ or $45312$).

                Conversely, suppose $m$ is the last iterate at which a
                crossing occurs in $D_{SW}(\sigma)$ (i.e.,
                $D_{SW}^{(\ell)}(\sigma)$ has no crossings for
                $\ell>m$).  We will prove that $L^{(m)}_2$ may have a
                crossing only at the first or last segment.  This, in
                turn, implies that row $m$ in $R$ or $S$ is
                decreasing.  A crossing occurs when there is a vertex
                of $L^{(m)}_1$ not in the shadow of any point of
                $L^{(m)}_2$.  We will prove that it can only be the
                first or last vertex.  Let $\{(s_1,r_1),
                (s_2,r_2),\dots\}$ and $\{(u_1,t_1),
                (u_2,t_2),\dots\}$ be the vertices that define
                $L^{(m)}_1$ and $L^{(m)}_2$, respectively.  Then
                $\{r_i\}_{i\ge 1}$ and $\{t_i\}_{i\ge 1}$ are
                decreasing while $\{s_i\}_{i\ge 1}$ and $\{u_i\}_{i\ge
                1}$ are increasing.  Write $(a,b)\le (c,d)$ if $(a,b)$
                is in the shadow of $(c,d)$ (i.e. if $a\le b$ and
                $c\le d$), and consider $L^{(m+1)}_1$ and
                $L^{(m+1)}_2$.  They are noncrossing and defined by
                points $\{(s_1,r_2), (s_2,r_3),\dots\}$ and
                $\{(u_1,t_2), (u_2,t_3),\dots\}$, respectively.  Then,
                for any $i$, $(s_i,r_{i+1})\le (u_j,t_{j+1})$ for some
                $j$.  Suppose $(s_i,r_{i+1})\le (u_j,t_{j+1})$ and
                $(s_{i+1},r_{i+2})\le (u_k,t_{k+1})$ for some $j<k$.
                Each upper shadowline vertex must contain some lower
                shadowline vertex in its shadow, so for all
                $\ell\in[j,k]$, $(s_i,r_{i+1})\le(u_\ell,t_{\ell+1})$
                or $(s_{i+1},r_{i+2})\le(u_\ell,t_{\ell+1})$.  Choose
                the least $\ell\in[j,k]$ such that
                $(s_{i+1},r_{i+2})\le(u_\ell,t_{\ell+1})$.  If
                $(s_i,r_{i+1})\le(u_\ell,t_{\ell+1})$, then
                $(s_{i+1},r_{i+1})\le (u_\ell,t_{\ell+1})\le
                (u_\ell,t_\ell)$.  If
                $(s_i,r_{i+1})\nleq(u_\ell,t_{\ell+1})$, then
                $(s_i,r_{i+1})\le(u_{\ell-1},t_\ell)$, so
                $(s_{i+1},r_{i+1})\le (u_\ell,t_\ell)$.  Thus, in both
                cases, $(s_{i+1},r_{i+1})\le (u_\ell,t_\ell)$, and the
                desired conclusion follows.
            \end{proof}\bigskip

            An immediate corollary of the above proof is that all rows
            $i\geq m$ in both $R(\sigma)$ and $S(\sigma)$ are monotone
            increasing from left to right if and only if every iterate
            $D_{SW}^{(i)}(\sigma)$ ($i\geq m$) is free from crossings.
            
            One can equivalently characterize intersecting shadowlines
            beyond the $0^{\rm th}$ iterate of $\sigma \in
            \mathfrak{S}_n$ in terms of sub-pile patterns for the
            entries in $R(\sigma)$ and $S(\sigma)$.  We state the
            following such result only for horizontal crossings, but
            vertical crossings can then be characterized by inverting
            $\sigma$ (i.e., by transposing within these pairs of
            patterns via a Sch\"{u}tzenberger-type symmetry result
            proven in \cite{refBLFPSAC05}).  Moreover, it is not
            difficult to show that avoiding both horizontal and
            vertical crossings in every iterate is equivalent to
            avoiding all crossings.

            \begin{cor}
                If $R(\sigma)$ and $S(\sigma)$ contain either of the
                following two simultaneous sub-pile patterns, then the
                permutation $\sigma \in \mathfrak{S}_{n}$ has a horizontal
                crossing in $D_{SW}^{(m)}(\sigma)$ (here $\{x_s\}_{s\ge 1}$ and
                $\{y_r\}_{r\ge 1}$ are monotone increasing; $m\le k,l$;
                and the numbers in the boxes indicate the number of
                elements in respective sub-piles):
            \end{cor}

                \hspace{-0.5in} \begin{tabular}{c c}

                    \begin{minipage}[c]{2.7in}

                        \begin{tabular}{l l}

                      \begin{tabular}{c c}
                          \fbox{$i$} & \\
                          $y_{1}$ & \fbox{$j$}\\
                          $y_{3}$ & $y_{2}$\\
                          \fbox{$k$} & \fbox{$m$} \\
                          \hline
                      \end{tabular}
                      $\subset R$
                      \raisebox{-1.125cm}{,}\!\!\!
                      &

                      \begin{tabular}{c c}
                          \fbox{$k-m$} & \\
                          $x_{1}$ & \fbox{$0$}\\
                          $x_{2}$ & $x_{3}$\\
                          \fbox{$i+m$} & \fbox{$j+m$} \\
                          \hline
                      \end{tabular}
                      $\subset S$

                      \end{tabular}

                    \end{minipage}

                    & \qquad or

                    \begin{minipage}[c]{2.7in}

                        \begin{tabular}{l l}

                      \begin{tabular}{c c}
                          \fbox{$i$} & \\
                          $y_{1}$ & \fbox{$j$}\\
                          $y_{3}$ & $y_{2}$\\
                          \fbox{$k$} & \fbox{$l$} \\
                          \hline
                      \end{tabular}
                      $\subset R$
                      \raisebox{-1.125cm}{,}\!\!\!
                      &

                      \begin{tabular}{c c}
                          \fbox{$k-m$} & \fbox{$l-m$}\\
                          $x_{2}$ & $x_{1}$\\
                          $x_{3}$ & $x_{4}$\\
                          \fbox{$i+m$} & \fbox{$j+m$} \\
                          \hline
                      \end{tabular}
                      $\subset S$

                      \end{tabular}

                    \end{minipage}

                  \end{tabular}

\end{document}